\documentclass[10pt]{article}

\title{ \Large \bf Generic  Hyperbolicity   for  the equilibria  \\ 
      of the one-dimensional parabolic equation 
    {\boldmath  $ u_t \: = \:  {(a (x) u_x)}_x + f(u).$} } 
 
\author{Antonio L. Pereira \thanks{ 
 Research partially supported by FAPESP-SP-Brazil, grant 2002/07464-8 }}

\date{ }

\newcommand{\Erre}{I \kern -3.0pt R} 
 
\newtheorem{lemma}{ Lemma} 
 \newtheorem{theo}{Theorem} 
\newtheorem{rem}{Remark} 
\newtheorem{cor}{Corollary} 
 \newtheorem{defin}{Definition}

\begin{document} 
\maketitle 
 
\vspace{-5mm} 
 
\begin{abstract} 
  We show, for some classes of diffusion coefficients that, 
 generically in $f$, all equilibria of the reaction-diffusion 
 equation 
 
   \[ 
    \begin{array}{lr} 
     u_t \: = \:  {(a (x) u_x)}_x + f(u)  &  0 < x  < 1 
     \end{array} 
  \] 
 with homogeneous Neumann boundary conditions are hyperbolic. 
 \end{abstract}

 \section{Introduction}

\noindent   We   consider here the  scalar parabolic equation with Neumann 
 boundary 
 conditions: 
 
 \begin{equation} \label{parabolic} 
 \left\{ 
 \begin{array}{lcr} 
   u_t \: = \:  {(a (x) u_x)}_x + f(u),  &  0 < x  < 1  &   \\ 
     u_x (0,t) \:= \: u_x (1,t) \: = \:0 &   & 
 \end{array} 
 \right. 
\end{equation} 
      where  $f: \Erre  \rightarrow \Erre $ is a $ {\cal C}^2 $ function and 
    $ a : [ 0,1 ]  \rightarrow {\Erre}^{+} = ( 0, \infty ) $ is 
    sufficiently smooth. 
 \vspace{2mm} 
 
  \noindent  An equilibrium point $u$ of (\ref{parabolic}) 
 is called {\em hyperbolic}  if and only if  0 is not one of the eigenvalues 
 of the linearized problem : 
 \[ 
    {(a w_x)}_x +  f'(u) w = \lambda w , \qquad w_x(0)=w_x(1) = 0. 
  \] 
\vspace{2mm} 
 
 \noindent  For {\em constant diffusion coefficient} $a$, it has been proved 
  by 
  Brunovsky and Chow (\cite{bc:hip}) that generically for $f$ in the 
  Whitney topology, all equilibria of  \ref{parabolic}  are hyperbolic. 
   The same result was obtained by Henry \cite{henry} under slightly more 
   general  hypotheses. 
\vspace{2mm} 
 
 \noindent For {\em nonconstant diffusion coefficient, } it has been shown by 
 Rocha 
  \cite{hr:nond}, that generically in the pair  $( a, f)$, all equilibria 
  are either hyperbolic or  non-degenerate ( with an appropriate definition 
  of the latter concept).  He also proved that the equilibria satisfying 
  an additional hypothesis can be made hyperbolic or non-degenerate, 
  by  perturbing 
  only $f$. We have shown  in \cite{pereira}  that generic hyperbolicity  can be be obtained {\em 
  with respect to $a$ }  and also with respect to $f$ {\em if 
 additional properties are imposed on the equilibria.}  Our purpose here is 
 to eliminate those additional requirements. To achieve this goal, we will 
 need to impose some extra hypotheses in the diffusion coefficient. 
 
  We  observe   that, if  explicit dependence of $f$ in $x$ 
 is allowed 
 such a result is  much easier  to obtain (see for example \cite{polacik}). 
    When the dependence on $x$ occurs only through $u$,  the result  is 
 certainly expected to hold 
  but to the 
 extent of our knowledge a proof is not yet available.

 We finally observe that, once hyperbolicity of equilibria has been obtained, 
  the Morse-Smale property follows  since transversality of the invariant manifolds is automatic 
 as proved in \cite{angenent} and 
  \cite{henry}.

 \section{Some technical  lemmas} 
 
 We state here some auxiliary results we will need in the proof.

    \begin{lemma}\label{sumzero} 
 
     Suppose $u:[a,b]\mapsto \Erre $ is a ${\cal C}^2$ function with nondegenerate 
 critical points 
    and  $\phi:[a,b]\mapsto \Erre$ 
  is continuous. . 
   Let 
 $R_q= \{ p\in [a,b] \ | \ u(p)=q, u'(p) \neq 0 \}$, 
 $C_q = \{ p\in [a,b] \ | \ u(p)=q, u'(p)= 0 \}$. 
  If  $ \int_a^b f(u) \phi \, dx =0 $, for any continuous   real function $f$. 
 then we have

\begin{equation} \label{eqnsumzero} 
 \sum_{p \in R_q \cap (a,b)} \frac{\phi(p)}{|u'(p)|} = 0. 
 \end{equation} 
and 
\begin{equation} 
 \sum_{p \in C_q \cap (a,b)}\frac {\phi(p)}{ \sqrt{|u''(p)|}} 
  + \frac{1}{2} \sum_{p \in C_q \cap \{a,b\}}\frac {\phi(p)}{ \sqrt{|u''(p)|}}  = 0. 
 \end{equation} 
 
 Reciprocally, if  (\ref{eqnsumzero}) holds for any $q \in R_q$ then 
   $ \int_a^b f(u) \phi \, dx =0 $, for any  continuous function  real function $f$. 
 \end{lemma}

  \noindent {\bf Proof. }  Let $p$ be a point where $u'(p) \neq 0$ and 
 $ q = u(p)$.  We can  find a local inverse 
 $u^{-1} : (q- \varepsilon, q+ \varepsilon) \to  (p - \delta_1, p + \delta_2) $, 
 where $\delta_1, \delta_2, \varepsilon $ are positive numbers.  We suppose 
 first that $u$ is increasing in  $(p - \delta_1, p + \delta_2) $ and, 
 therefore, $ u(p- \delta_1) = q - \varepsilon,  u(p = \delta_2) = q + \varepsilon $. 
 
 Let   $f_{\varepsilon}  : \Erre \to \Erre $ be a ${\cal C}^{\infty}$ function 
 satisfying 
\[ f(x) = \left\{ 
\begin{array}{cc} 
 \frac{1}{\varepsilon} & \mbox{ if  }  | y - q | \leq \varepsilon \\ 
   \leq \frac{1}{\varepsilon}  & \mbox{ if  } 
    \varepsilon \geq  | y - q | \leq  \varepsilon  + \varepsilon^2 \\ 
   0 & \mbox{ otherwise. } 
     \end{array} \right. 
  \] 
 
 Then we have 
 
   \begin{eqnarray*} 
 \int_{p- \delta_1}^{p + \delta_2} f_{\varepsilon}(u(x)) \phi(x) \, d \, x 
  & = & 
  \int_{q- \varepsilon}^{q + \varepsilon} f_{\varepsilon}(y) \phi(u^{-1}(y)) 
   \cdot \frac{1}{u'(u^{-1}(y))} \, d \, y \\ 
  & = & 
   \int_{q- \frac{\varepsilon}{2}}^{q + \frac{\varepsilon}{2}} \frac{1}{\varepsilon} 
  \phi(u^{-1}(y)) \frac{1}{u'(u^{-1}(y))}\, d \, y  \\ 
  &  +  & 
   \int_{q-\frac{\varepsilon}{2} - {\varepsilon^2}}^{q-\frac{\varepsilon}{2}} 
  f_{\varepsilon}(y)  \phi(u^{-1}(y)) \frac{1}{u'(u^{-1}(y))} \, d \, y \\ 
  &  + & 
\int_{q+\frac{\varepsilon}{2}}+{q+\frac{\varepsilon}{2} + {\varepsilon^2}} 
  f_{\varepsilon}(y)  \phi(u^{-1}(y))  \frac{1}{u'(u^{-1}(y))} \, d \, y   \\ 
     \end{eqnarray*}

  It easy to see  that the last two integrals go to zero as $ \varepsilon 
 \to 0$. For the first, we have

 \begin{eqnarray*} 
   \int_{q- \frac{\varepsilon}{2}}^{q + \frac{\varepsilon}{2}} \frac{1}{\varepsilon} 
  \phi(u^{-1}(y)) \frac{1}{u'(u^{-1}(y))} \, d \, y   & = & 
 \frac{1}{\varepsilon} \cdot \varepsilon \cdot  \phi(u^{-1}(\xi))  {u'(u^{-1}(\xi))} 
    \end{eqnarray*} 
 where $\xi \in [ {q- \frac{\varepsilon}{2}}, {q + \frac{\varepsilon}{2}} ].$ 
 
  As $ \varepsilon \to 0$ we then obtain 
 
 \begin{eqnarray*} 
 \int_{p- \delta_1}^{p + \delta_2} f_{\varepsilon}(u(x)) \phi(x) \, d \, x 
 &  \to &  \phi(u^{-1}(q)) \frac{1}{u'(u^{-1}(q))} \\ 
  & = & \frac{ \phi(p) }{u'(p)} \\ 
  & = & \frac{ \phi(p) }{ | u'(p) |}. 
 \end{eqnarray*}

  If $u$ is decreasing  in  $(p - \delta_1, p + \delta_2) $ we obtain 
  similarly that 
 
   \begin{eqnarray*} 
 \int_{p- \delta_1}^{p + \delta_2} f_{\varepsilon}(u(x)) \phi(x) \, d \, x 
 &  \to &  -  \phi(u^{-1}(q))  \frac{1}{u'(u^{-1}(q))} \\ 
  & = &  - \frac{ \phi(p) }{u'(p)} \\ 
  & = & \frac{ \phi(p) }{ | u'(p) |}. 
 \end{eqnarray*}

   Now, let  $q$ be a regular value of $u$ and  denote by $I_q$ the pre-image 
 of $q$ by $u$, that is 
 $I_q  = 
 \left\{ p \in [a,b] \  |  \ u(p) = q   \right\} $. 
 
     Then, with $f_{\varepsilon}$ as above, we obtain

\begin{eqnarray*} 
0 & = & \int_{a}^{1} f_{\varepsilon}(u(x)) \phi(x) \, d \, x \\ 
  & = & 
 \sum_{p \in I_q} 
 \int_{p- \delta_1}^{p + \delta_2} f_{\varepsilon}(u(x)) \phi(x) \, d \, x  \\ 
  & \to & 
  \sum_{p \in I_q} \frac{ \phi(p) }{ | u'(p) |}. 
  \end{eqnarray*} 
 (Of course the values of  $\varepsilon$ above depend on the $p$. We chose not 
 to write this dependence explicitly to simplify  the      notation).

 Suppose now that $p\in (a,b)$ is a critical point of $u$. 
  We can write $u$  in a neighborhood of $p$  as 
  \begin{eqnarray*} 
   u(x) & = & u(p) + \frac{u''(p)}{2} (x-p)^2 + O(x-p)^3. 
 \end{eqnarray*} 
 
 Observe that $u''(p) \neq 0$ by hypothesis. 
 Suppose  $u''(p)  < 0$.  The image of a (small) interval $I$ around 
 $p$ is then an interval of the form $ [q - \varepsilon, q] $, where 
 $q = u(p).$  We have, for $x \in I$, $ y  \in  [q - \varepsilon, q] $ 
 
  Write 
   $\alpha = \frac{u''(p)}{2}, w = (x-p)^2, z = y-q.$  Then

 \begin{eqnarray*} 
 u(x)  =  y  &  \Leftrightarrow &  \alpha w + O(|w|^{3/2})  = z  \\ 
 &  \Leftrightarrow &  w = \frac{z}{ \alpha + O(|w|^{1/2})} \\ 
 & = & \frac{z}{ \alpha + O(|z|^{1/2})} \\ 
 & = & \frac{z}{\alpha} \frac{1}{ 1 + O(|z|^{1/2})} \\ 
 & = &  \frac{z}{\alpha}  \left( { 1 + O(|z|^{1/2})} \right) \\ 
 & = &  \frac{z}{\alpha}   + O(|z|^{3/2}) 
 \end{eqnarray*} 
 
 Therefore $ u$ is  one to one from each one of the intervals 
$[ p- \delta_1, p]$ and 
 $[p, p + \delta_2]$  into $ [q - \varepsilon, q] $, with

\begin{eqnarray*} 
 (x-p) & = & \pm \sqrt{ \frac{y-q}{\alpha} + O( |y-q|^{3/2}) } \\ 
       & = & \pm \sqrt{ \frac{y-q}{\alpha}} 
              \sqrt{1 + O( |y-q|^{1/2}) } \\ 
      & = & \pm \sqrt{ \frac{y-q}{\alpha}} 
              \left(1 + O( |y-q|^{1/2}) \right) \\ 
    & = &  \pm \sqrt{ \frac{y-q}{\alpha}} 
       \left(   1 + O( |y-q|^{1/2})  \right)  \\ 
     & = & \pm \sqrt{ \frac{y-q}{\alpha}} + 
               + O( |y-q|) . 
 \end{eqnarray*}

 Let   $g_{\varepsilon}  : \Erre \to \Erre  $ be a ${\cal C}^{\infty}$ function 
 satisfying 
\[ g_{\varepsilon} (x) = \left\{ 
\begin{array}{cc} 
 \frac{1}{\sqrt{\varepsilon}} & \mbox{ if  }  | y - q | \leq \varepsilon \\ 
   \leq \frac{1}{\sqrt{\varepsilon}}  & \mbox{ if  } 
    \frac{\varepsilon}{2} \geq  | y - q | \leq 
     \frac{ \varepsilon  + \varepsilon^2}{2} \\ 
   0 & \mbox{ otherwise. } 
     \end{array} \right. 
  \] 
 
 Denoting by $x(y)$ the inverse of $u: [ p- \delta_1, p]  \to 
 [q - \varepsilon, q].  $, we obtain

   \begin{eqnarray*} 
 \int_{p- \delta_1}^{p } g_{\varepsilon}(u(x)) \phi(x) \, d \, x 
  & = & 
  \int_{x(q- \varepsilon/2 - {\varepsilon^2}/2)}^{x(q- \varepsilon/2 ) } 
  g_{\varepsilon}(u(x)) \phi(x)   \, d \, x + 
  \int_{x(q- \varepsilon/2 )}^{p } g_{\varepsilon}(u(x)) \phi(x) 
   \, d \, x \\ 
    \end{eqnarray*}

   For the first integral, we  have 
 
    \begin{eqnarray*} 
  \left\| \int_{x(q- \varepsilon/2 - {\varepsilon^2}/2)}^{x(q- \varepsilon/2 ) } 
  f_{\varepsilon}(u(x)) \phi(x)   \, d \, x  \right\| 
   &  \leq &  \frac{1}{\sqrt{\varepsilon}} ||\phi||_{\infty} 
    \left( \sqrt{\frac{-\varepsilon }{2\alpha }} - 
          \sqrt{\frac{-\varepsilon -\varepsilon^2 }{ 2\alpha }} + O(\varepsilon)  \right) 
  \\ 
  & = &  \frac{1}{\sqrt{\varepsilon}} O(\varepsilon) \\ 
  & = &  O(\sqrt{\varepsilon}) 
   \end{eqnarray*}

 For the second integral, we have 
 
  \begin{eqnarray*} 
   \int_{x(q- \varepsilon/2 )}^{p }g_{\varepsilon}(u(x)) \phi(x)  \, d \, x 
    & = &  \frac{1}{\sqrt{\varepsilon}}(p - x(q-\varepsilon/2) )  \phi(\bar{x}) \\ 
   & = &  \frac{1}{\sqrt{\varepsilon}}( \sqrt{\frac{-\varepsilon}{2\alpha} } ) 
   \phi(\bar{x})  \\ 
  & = &  \frac{ \phi(\bar{x}) }{ \sqrt{-2\alpha} }. 
   \end{eqnarray*} 
 where $\bar{x} \in [q-\varepsilon/2, q]$.

   As $ \varepsilon \to 0$ we  obtain 
 
\begin{eqnarray*} 
  \int_{p- \delta_1}^{p } f_{\varepsilon}(u(x)) \phi(x) \, d \, x 
&  \to &  \frac{ \phi(p) }{ \sqrt{-2\alpha} }. 
   \end{eqnarray*}

 Similarly, we obtain

\begin{eqnarray*} 
  \int^{p+ \delta_2}_{p }g_{\varepsilon}(u(x)) \phi(x) \, d \, x 
&  \to &  \frac{ \phi(p) }{ \sqrt{-2\alpha} }. 
   \end{eqnarray*}

Therefore, 
 
  \begin{eqnarray*} 
  \int^{p+ \delta_2}_{p-\delta_1 }g_{\varepsilon}(u(x)) \phi(x) \, d \, x 
&  \to & \sqrt{2} \frac{ \phi(p)}{ \sqrt{|\alpha|} }. 
   \end{eqnarray*} 
 
as $ \varepsilon \to 0$, and the same is true if  $\alpha = \frac{u''(p)}{2} >0$.

 If $p = a$ or $p= b$, the same computations can be done  in the intervals 
 to the right or left respectively, giving the same result divided by $2$

 Now, let  $q$ be a critical  value of $u$ and  denote by $I_q$ the pre-image 
 of $q$ by $u$, that is 
 $I_q   = 
 \left\{ p \in [a,b] \  |  \ u(p) = q  \right\} $. 
 
 $ C _q $ the critical points in $I_q$ and   $R_q $ the regular points in $I_q$, 
 that is 
 
 $C_q  = 
 \left\{ p \in [a,b] \  |  \ u(p) = q, u'(p) = 0  \right\} $

and   $R_q = I_q - C_q $ the regular points in $I_q$,

     Then, with $g_{\varepsilon}$ as above we obtain

\begin{eqnarray*} 
0 & = & \int_{a}^{b}g_{\varepsilon}(u(x)) \phi(x) \, d \, x \\ 
  & = & 
 \sum_{p \in I_q} 
 \int_{p- \delta_1}^{p + \delta_2} g_{\varepsilon}(u(x)) \phi(x) \, d \, x  \\ 
 & = & 
 \sum_{p \in C_q} 
 \int_{p- \delta_1}^{p + \delta_2} g_{\varepsilon}(u(x)) \phi(x) \, d \, x + 
 \sum_{p \in R_q} 
 \int_{p- \delta_1}^{p + \delta_2} g_{\varepsilon}(u(x)) \phi(x) \, d \, x   \\ 
 & = & 
 \sum_{p \in C_q} 
 \int_{p- \delta_1}^{p + \delta_2} g_{\varepsilon}(u(x)) \phi(x) \, d \, x. 
  \end{eqnarray*} 
 (Of course the values of  $\delta_1,\delta_2 $ above depend on  $p$. 
 If $p=a$ (resp. $p=b$), then  $\delta_1 = 0 $ resp. $\delta_2=0$ 
  We chose not 
 to write this dependence explicitly to simplify the  notation).

 As $ \varepsilon \to 0$ we obtain

\begin{eqnarray*} 
  \sum_{p \in C_q \cap (a,b)} \sqrt{2} \frac{ \phi(p)}{ \sqrt{|u''(p)|}} 
 +\sum_{p \in C_q \cap \{a,b\}} \frac{\sqrt{2}}{2} 
  \frac{ \phi(p)}{ \sqrt{|u''(p)| }} = 0. 
 \end{eqnarray*}

 To prove the converse, we first denote by  $ 0=x_0 < x_1 < \cdots < x_n=1$ the critical 
 points of $u$ in $ [0,1]$, and by $ y_0 < y_1 < \cdots < y_m$ the critical values. 
 
 Let $I_j = [y_{j-1}, y_{j}]. $For a fixed $j$  and any  $y$ in the interior 
 of $I_j$ let 
   $n_{i_j}$ be the number  of points in $u^{-1}(y)$ (observe that  $u^{-1}(I_j)$ is then 
 the disjoint union of $n_{i_j}$ intervals). 
 By hypothesis, we have 
 
 \[    \sum_{i=1}^{n_{i_j}}   \frac{ \phi(\xi_i(y))}{  |u'( \xi_i(y))  | }    = 
  \sum_{i=1}^{n_{i_j}}   \phi(\xi_i(y)) \cdot | \xi_i(y)  | = 0. \] 
 
 To simplify the notation, we write this sum as 
 
  \[ 
  \sum_{i=1}^{n}   \phi(\xi_i(y)) \cdot | \xi_i(y)  | = 0. \] 
 with the understanding that  the summand is zero, if  $ y $  is not in the range of $\xi$. 
 
 If $g$ is a continuous  function in $I_j$,  we then have 
 
 \begin{eqnarray} 
 0 & =  &    \sum_{i=1}^{n}  \int_{I_j} g(y) \phi(\xi_i(y)) \cdot | \xi_i(y)  | \, dy \nonumber  \\ 
   &   =  &    \sum_{i=1}^{n}  \int_{\xi(I_j)} g(u(x)) \phi(x) \, dx 
 \end{eqnarray}

 Now, if $I_i = [x_{i-1}, x_i], \  i=1,2, \cdots, n$, we have 
 $I_i =\stackrel{\circ}{ \bigcup} \xi(I_j)$. Therefore, 
 
\begin{eqnarray} 
 0 & =  & \sum_{j=1}^m   \sum_{i=1}^{n}  \int_{\xi(I_j)} g(u(x)) \phi(x) \, dx  \nonumber \\ 
   &  =  &   \sum_{i=1}^{n}  \int_{(I_i)} g(u(x)) \phi(x) \, dx  \nonumber \\ 
   & = &   \int_{0}¹ g(u(x)) \phi(x) \, dx. 
 \end{eqnarray}

 \vspace{3mm}

  If $u$ is a (fixed) function in $[0,1]   $  and $ [a,b] \subset [0,1]   $  , we will say that 
 $\Phi: [0,1] \mapsto \Erre $ is in the 
 {\em orbit} of $u$ in  $ [a,b]    $, and   write $ \Phi \in {\cal O}(u)\{a,b\}  $  if 
  $\Phi(x) = g(u(x))$, for some continuous function $g:  \Erre \mapsto \Erre$ 
 for any $x \in [a,b]$. We will  also say that $ \Phi$ belongs to the  {\em q-orbit} of $u$ in 
 $ [a,b]$ if 
 there exists a function $g:  \Erre \mapsto \Erre$ 
 continuous, except maybe at the  critical points of $u$ with  $\Phi(x) = g(u(x))$ 
 for $x$ regular point of $u$. In this case, we  write 
 $ \Phi \in q-{\cal O}(u)\{a,b\}  $

  \begin{lemma} \label{orbit} Suppose $u$  is    analytic   in the interval $[a,b]$, with 
 nondegenerate critical points,   $u(a)  = u(b), \  u'(a) =0 $ (or $  u'(b) = 0$)  and, 
     if $p_1$, $p_2$ are   critical points of $u$ in $(a,b)$, 
 then  $u(p_1) \neq u(a)$  (or $u(p_1) \neq u(b) $)  and   $u(p_1) \neq u(p_2))$. 
  Suppose also that  $ \Phi$ is analytic 
 in the interval $[a,b]$  except maybe at the critical points of $u$. Then, 
  with the 
 notation of lemma (\ref{sumzero}), 
 if  $ \sum_{p \in R_q \cap (a,b)} \Phi(p) \frac{u'(p)}{|u'(p)|} = 0$, 
 $ \Phi$ belongs to the  {\em q-orbit} of $u$ in  $ [a,b]$.  In particular, if 
  $ \Phi$ is analytic  in the interval $[a,b]$  then  $ \Phi \in {\cal O}(u)\{a,b\}  $. 
 \end{lemma} 
 
 \noindent {\bf Proof. }   Let $ a = x_0 < x_1 <  x_2 < \cdots < x_n = b$ be the critical points 
 of $ u$ in $[a,b]$, and $\xi_i$, $ i=1,2 \cdots, n$ be the inverse of 
 $u$ in the interval $[x_{i-1},x_{i} ]$.  If $q$ is a regular value of $u$ and 
 $ \xi_{i_1},  \xi_{i_2}, \cdots, \xi_{i_n}$  are its inverse images under $u$, we 
 have, by hypotheses 
 $ \sum_{j=1}^{j=k} \Phi(\xi_{i_j}) \, \mbox{\rm sign} {u'(\xi_{i_j})} = 0$. 
 Since the critical points are nondegenerate they must be  local minima and maxima 
 in succession. 
  Suppose, for definiteness, that $x_0, x_2, \cdots, x_n$ are local  maxima and 
 $x_1,  \cdots, x_{n-1}$ are local minima. 
 Let $ i_1, i_2, \cdots, i_k$ 
  $k ={ (n-1)  \over 2}$ be such  that 
  $u(x_{i_1}) < u(x_{i_2}) < , \cdots, < x_{i_k}$. 
  Now, for $ u(x_{i_1}) < q < u(x_{i_2})$, we must have, by hypotheses, 
  $ - \Phi \circ \xi_{i_1}(q) +  \Phi \circ \xi_{i_1+ 1 }(q) = 0 $, that is 
 $  \Phi \circ \xi_{i_1}(q) =  \Phi \circ \xi_{i_1+ 1 }(q) $. 
  Since  $  \Phi \circ \xi_{i_1} $ and $  \Phi \circ \xi_{i_1+ 1 } $ are 
 analytic functions in their common (open) interval of definition, we must have  the equality 
 in this whole interval. Proceeding upward, we next find that 
  $ - \Phi \circ \xi_{i_1}(q) +  \Phi \circ \xi_{i_1+ 1 }(q) 
 - \Phi \circ \xi_{i_2}(q) + 
  \Phi \circ \xi_{i_2+ 1 }(q) = 0 $. Since we already know that 
 $  \Phi \circ \xi_{i_1}(q) =  \Phi \circ \xi_{i_1+ 1 }(q) $ this implies 
 $  \Phi \circ \xi_{i_2}(q) =   \Phi \circ \xi_{i_2+ 1 }(q) $. In this way we find 
 successively  that 
   $  \Phi \circ \xi_{i_1} =  \Phi \circ \xi_{i_1+ 1 },   \Phi \circ \xi_{i_2} = 
 \Phi \circ \xi_{i_2+ 1 }, \cdots, 
    \Phi \circ \xi_{i_{k-1}} =  \Phi \circ\xi_{i_{k-1}+1 }  $ or, returning 
 to the old ordering 
   $  \Phi \circ \xi_{1} =  \Phi \circ \xi_{2 },   \Phi \circ \xi_{3} = 
 \Phi \circ \xi_{4 }, \cdots, 
    \Phi \circ \xi_{n-1} =  \Phi \circ \xi_{n }$.

  Now, starting from the maximum value and going downward, we find similarly that 
  $  \Phi \circ \xi_{2} =  \Phi \circ \xi_{3 },   \Phi \circ \xi_{4} = 
  \Phi \circ \xi_{5 }, \cdots, 
    \Phi \circ \xi_{n-2} =  \Phi \circ \xi_{n-1 }$ and also that $  \Phi \circ \xi_{1} = 
 \Phi \circ \xi_{n }   $.  From this, it follows  that 
   $ \Phi \circ \xi_{i}(q) =  \Phi \circ \xi_{j}(q) $ for any $ 1 \leq i,j \leq n $ and any $q$ 
 that belongs  to the common interval of 
  $ \Phi \circ \xi_{i} $ and  $ \Phi \circ \xi_{j}$.  In other words, we must have 
  $ \Phi(p_1) = \Phi(p_2)$ whenever 
  $u(p_1) = u(p_2)$ and $ p_1, p_2$ are not critical points, which proves the claim.



  \begin{lemma} \label{inorbit} Suppose $a$ and $f$ are analytic functions,  $u$ is 
 a non-constant solution of (\ref{bvp})  and $a,b$ are points in $ [0,1]$ such that 
   $u(a)  = u(b), \  u'(a) =0  $ (or $ u'(b) = 0$)  and, 
     if $p_1$, $p_2$ are   critical point of $u$ in $(a,b)$, 
 then  $u(p_1) \neq u(a)$ (or $u(p_1) \neq u(b)  $) and   $u(p_1) \neq u(p_2))$ 
 Suppose there exists  a solution    $\eta$ 
   of the second order O.D.E. in (\ref{lve}), with 
  $\int_a^b g(u) \eta \, dx = 0 $, for any continuous real function $g$. 
 Then,  the following functions belong to  $ {\cal O}(u) \{a,b\} $ 
 \begin{eqnarray} &  &        a \eta_x   \label{02}  \\ 
                  &  &          au_x \eta  \label{03}  \\ 
                  &  &            a \eta^2  \label{04}  \\ 
                  &  &               a {u_x}^2 \label{05}  \\ 
                  &  &                a_x \eta \label{06}  \\ 
                  &  &                 a_x u_x \label{07} \\ 
                  &  &                   a u_{xx} \label{08} \\ 
                  &  &               a_{xx}   \label{11} \\ 
                   &  &               a   \label{12} \\ 
                   &  &                \eta_x   \label{14} \\ 
                    &  &               {\eta}^2   \label{15} \\ 
                     &  &               {u_x}^2   \label{16} \\ 
                      &  &               u_{xx}   \label{17} 
    \end{eqnarray} 
 Also the following functions belong to  $q-{\cal O}(u) \{a,b\} $ 
 
 \begin{eqnarray} &  &   \frac{\eta}{u_x} \label{01} \\ 
                  &  &   \frac{ \left( a_{xx} \eta + a_x \eta_x 
                          \right) }{u_x}   \label{09}  \\ 
                  &  &   a_{xx} +  \frac{a_x \eta_x }{\eta}     \label{10} \\ 
                  &  &   \frac{a_x  }{u_x}     \label{13} 
                    \end{eqnarray} 
 
 \end{lemma} 
 
 \noindent {\bf Proof. }  We suppose that $u'(a) = 0$. The case $u'(b) = 0$ is similar. 
  First, observe that (\ref{01})  follows immediately from  lemmas \ref{sumzero} 
and  \ref{orbit}.

   Multiplying (\ref{lve}) by any  ${\cal C}^1$ function    $   g(u)  \in   {\cal O}(u)$ 
 and integrating by parts, we obtain 
 
 \begin{eqnarray*} 
 0 & = &  \int_a^b  g(u) (a \eta_x)_x  u_x  \, dx  \\ 
   & = &  g(u)  a \eta_x |_a^1 -  \int_a^1  g'(u) a \eta_x  u_x  \, dx \\ 
 \end{eqnarray*} 
 Since $u(a) = u(b)$, for an arbitrary continuous function $h$, we can find a 
  $({\cal C}^1)$ 
 function $g$, with 
  $g(u(a)) = g(u(b)) = 0 $, and $g'= h$.

 Therefore, lemma \ref{sumzero}   applies to  $ a \eta_x  u_x$  in the place of 
 $ \phi$ and 
  the result follows 
 from lemma (\ref{orbit}). 
 
 Also, taking now $g(u(a)) \neq 0$, we obtain 
    $ a \eta_x  (a) =   a \eta_x  (b) $.

  Multiplying (\ref{bvp}) by any  ${\cal C}^1$ function 
   $   g(u)  \in   {\cal O}(u)$  times $\eta$, with $g(u(a)) = 0 $ 
  and integrating by parts, we obtain 
 
 \begin{eqnarray*} 
 0 & = & \int_a^1 \left(au_x \right)_x   g(u)  \eta +  f(u) g(u) \eta   \, dx \\ 
   & = &  -  \int_a^1  g'(u) \left( au_x  \right) u_x  \eta +  \left( a\eta_x  \right) g(u) u_x   \, dx \\ 
   & = &  -  \int_a^1  g'(u) \left( a u_x \eta  \right) u_x 
 \end{eqnarray*} 
 where, to obtain the last equality, we used that  $   a \eta_x \in {\cal O}(u)$. 
   Therefore (\ref{03}) also follows  from  lemmas \ref{sumzero} and 
   (\ref{orbit})  that 
  $ a u_x  \eta \in {\cal O}(u) $. 
 To obtain  (\ref{05}), we multiply  $ a u_x  \eta $ by 
 $\frac{u_x}{\eta}$,  and   (\ref{03})  and  (\ref{01}) to obtain the result at the 
 points where $ \eta \neq 0$, and then conclude $ u_x^2 \in {\cal O}{u} $ by continuity. 
 Note that, in particular, we must have $u_x(b) = 0$. 
 
 To obtain (\ref{06}) we can proceed in various ways.  We can, for instance,   multiply  (\ref{02}) 
 by $g(u) \eta$, integrate by parts and use (\ref{03}) (\ref{04}) to obtain 
 $ \frac{a_x \eta^2}{u_x} \in   {\cal O}(u) $ and then use  (\ref{01}). 
 
  Now  (\ref{07}) and  (\ref{08}), follow from  (\ref{06}),  (\ref{01}) and  equation (\ref{bvp}) 
 
  (\ref{09}) can also be obtained in various ways. For example, from (\ref{06}) 
 
 $ a_x \eta \in    {\cal O}(u)  \Rightarrow g(u) f'(u)  a_x \eta \in    {\cal O}(u) 
  \Rightarrow g(u)   a_x  \left(a \eta_x  \right)_x \in    {\cal O}(u)$ for any $g(u)$. 
  Multiplying by $\eta$  we arrive after integrating by parts and 
  using    (\ref{06}) at 
   $ \int_0^1   \left(  g(u) a \eta_x \right) \left( a_{xx} \eta + a_x \eta_x 
   \right)   \, dx= 0 $. 
 
   Since  $ g(u) a \eta_x$ is arbitrary, except at the points where $\eta_x = 0$ the conclusion of lemma \ref{sumzero} 
 (with  $ \left( a_{xx} \eta + a_x \eta_x  \right)$ in the place of $\phi$)   can still be obtained 
  at least at these  points. But, it  must also hold  at all points where $ u_x \neq 0 $  so  (\ref{orbit}) can be applied 
 again to obtain the result. 
 
  Multiplying  (\ref{09}) by $\frac{u_x}{\phi}  $, we obtain   (\ref{10}). 
 
  Dividing   (\ref{02}) by  (\ref{04}) and then multiplying by    (\ref{06})  we obtain that 
  $  \frac{a_x \eta_x }{\eta}   \in q-{\cal O}(u) \{a,b\}  $ and so 
   (\ref{11}) follows from   (\ref{10}). 
 
  Now multiplying $a_{xx}$ by $ g(u) u_x$, it follows that 
  \begin{eqnarray*} 
 0 & = & \int_a^b g(u) a_{xx} u_x   \, dx \\ 
   & = &  -  \int_a^b  g'(u) u_x^2 a_x   \, dx -  \int_a^b  g(u) u_{xx} a_x    \\ 
   & = &  -   \int_a^b  g(u) u_{xx} a_x 
 \end{eqnarray*} 
 
Therefore $ \frac{a_x u_{xx}}{u_x}  \in q-{\cal O}(u) \{a,b\}  $. From this, and  (\ref{08}), we conclude that there is a 
 function $g(u)$, with $ \frac{a_x}{a}  = g(u) u_x  $, at least at the points 
 where $u_x \neq 0$. 
 Now, we claim that  $a_x= 0$ if $u_x=0$. In fact, from lemma 
 \ref{sumzero}, we know that $n=0$ at those points and, from (\ref{06}), 
 $a_{xx} \eta + a_x \eta_x = 0$. If $a_x \neq 0$, it follows then that 
 $ \eta(x) = \eta(x) = 0$ and $ \eta \equiv 0$ a contradiction. 
  Therefore  the equality must hold everywhere. 
 Integrating, we obtain 
 $ \ln(a) = G(u)$ where $G$ is a primitive of $g$, from which (\ref{12}) and   (\ref{13}) follow immediately. 
  (\ref{14}), (\ref{15}), (\ref{16})  and (\ref{17}) follow then from  (\ref{12}) 
 and  (\ref{02}), (\ref{04}), (\ref{05}) and  (\ref{08})  respectively.

 \section{Hyperbolicity of the equilibria}

 \noindent Let  $f: \Erre  \rightarrow \Erre $ be a $ {\cal C}^2 $ function and 
    $ a : [ 0,1 ]  \rightarrow {\Erre}^{+} = ( 0, \infty ) $ 
    continuous. 
 
 \noindent  We denote by  $ {  E(a,f)}$  the set of equilibria of 
 (\ref{parabolic}). Clearly,  $u \in {  E(a,f)} $ if and only if it is 
   a solution of 
   the boundary value 
 problem

  \begin{equation} 
  \left\{ 
    \begin{array}{lr} 
      {\displaystyle {(a (x) u_x)}_x + f(u) \: = 0 }, \: &  0 \leq x \leq 1 \\ 
     {\displaystyle u_x (0)  =   u_x (1) = 0} & \ 
     \end{array} 
   \right.   \label{bvp} 
    \end{equation}

 The  initial value problem associated to (\ref{bvp}) is

 \begin{equation} 
  \left\{ 
    \begin{array}{lr} 
      {\displaystyle {(a (x) u_x)}_x + f(u) \: = 0 }, \: &  0 \leq x \leq 1 \\ 
     {\displaystyle u_x (0)=0  \quad  u (0) =  u_0 } & \ 
     \end{array} 
   \right.   \label{ivp} 
    \end{equation} 
 and 
  the linear variational equation corresponding to 
 (\ref{bvp})  around  $u$ is 
 
  \begin{equation} 
  \left\{ 
    \begin{array}{lr} 
      {\displaystyle {(a (x) \phi_x)}_x + f'(u)\cdot \phi 
       \: = 0 }, \: &  0 \leq x \leq 1 \\ 
     {\displaystyle \phi_x (0) \:= \: \phi_x (1) \: = \:0} & \ 
     \end{array} 
   \right.   \label{lve} 
    \end{equation}

  A solution of (\ref{bvp}) is therefore hyperbolic if, and only if, 
 (\ref{lve}) has only the trivial solution.

\noindent It is easy to show (see \cite{hr:nond})  that,  if the set of 
 solutions of 
 (\ref{bvp}) is non-empty, then it contains a constant solution $\xi$. Also, 
 since the change of variables $ u ' = u- \xi$ does not affect the 
 boundary conditions on (\ref{parabolic}) we can, and will, 
 assume that $f(0)=0$.

 In what follows we denote by ${\cal C}^2_S ( \Erre ) $ the space of 
   ${\cal C}^2$ functions $ f: \Erre \rightarrow \Erre$  satisfying 
 $f(0) = 0$ with the  Whitney 
 topology.

  Let $ H_N^2[0,1] = \{ u \in H^2[0,1] | | \ u_x(0) = u_x(1) = 0    \} $ and 
 consider the map

  \begin{eqnarray*} 
\Psi_f: {H}^2_N [0,1] & \mapsto & \mathnormal{L}^2[0,1] \\ 
{} u & \mapsto  & (au_x)_x+f(u) 
\end{eqnarray*}

 \begin{lemma}\label{regular} 
  An equilibrium of  (\ref{parabolic}) is  hyperbolic if, and only if, 
 it is   a regular  point of $\Psi_f$. 
 \end{lemma} 
 \noindent {\bf Proof. }  $u$ is a regular point of $\Psi_f$ if and only if, 
  the derivative 
 \[ 
  D_u \Psi_f {H}^2_N [0,1]  \mapsto  {L}^2[0,1] 
: \phi \longrightarrow (a\phi_x)_x+f'(u)\phi 
 \] 
  is surjective. Now  $D_u \Psi_f$ is Fredholm of index $0$ and, therefore, 
 it is surjective if and only it is an isomorphism. 
 
 \vspace{3mm}

 Let us consider the equilibria of  (\ref{parabolic}) with 
 norm less than or equal to $N$, for some natural $N$. Our plan is to prove 
 that they are all hyperbolic for an open dense set of $f$ in the 
 Whitney topology and then take intersection. Now, for these equilibria, 
 we can use the ${\cal C}^2$ topology for $f$ since its  values 
 outside a compact set are irrelevant.

 We  first  look at the {\em constant equilibria}. 
 
  \begin{lemma} \label{constequilibria} 
 The constant equilibria $u$ of (\ref{parabolic}), with $ ||u||_{H^2} \leq N$ 
 are all hyperbolic for $f$ in a open dense set of ${\cal C}^2$. 
\end{lemma} 
 
 \noindent {\bf Proof. }  As usual,  openness    is not a problem here, so we just 
 have to prove density. Now, a constant $u_0$ is a non hyperbolic  equilibrium 
 of (\ref{parabolic}) if and only  $f(u_0) =0$ and $f'(u_0)$ is one of 
 the eigenvalues of the operator 
  \[ \phi \mapsto \left( a \phi_x  \right)_x. \] 
  in $L^2[0,1]$. We can first choose $f$ with only a finite number of 
 zeroes in  the  ball of radius $N$ and center at  the origin, 
 so that the number of constant equilibria 
  with $ ||u||_{H^2} \leq N$ is finite. 
  Now, since the eigenvalues of the above operator form a discrete set it 
 is easily seen that $f'$ can be modified, without changing the zeroes of $f$, 
 in such a way that $f'(u_0)$ will not be an eigenvalue whenever $u_0$ is 
 a zero.

   Let now $Y_N$ be the open dense  dense set in ${\cal C}^2$, given by 
 lemma (\ref{constequilibria}), 
 $ B_N = \{ u \in {H}^2_N [0,1] \ | \ ; \|u\|_{{H}^2_N} [0,1]  \leq N $, 
 and consider the map 
 \begin{eqnarray}\label{Psi} 
 \Psi:B_N \times Y_N & \mapsto & {L}^2[0,1] \\ 
(u,f) & \mapsto  & (a u_x)_x + f(u) 
\end{eqnarray}

 \begin{lemma}\label{condition} 
  Suppose $0$ is a not a regular value of the $\Psi$. Then there is an 
 equilibrium $u$ of (\ref{parabolic}) and a  nontrivial solution $\phi$ of the corresponding 
 linearized equation   (\ref{lve}) such that 
 $\int_0^1 \dot{f}(u)\phi = 0 \textrm{ } \forall \dot{f} \in 
  {\cal C}^2(\Erre, \Erre)$. 
 \end{lemma} 
 
 \noindent {\bf Proof. } Suppose there exists 
  $(u,f) \in F^{-1}(0)$  such that 
\begin{eqnarray*} 
D\Psi(u,f): {H}^2_N[0,1] \times {\cal C}^2(\Erre,\Erre) & 
\mapsto & {L}^2[0,1] \\ 
(\dot{u},\dot{f}) & \mapsto & (a \dot{u}_x)_x + f'(u)\dot{u} 
+\dot{f}(u) 
\end{eqnarray*} 
is not surjective.

 By hypothesis, there exists    a nontrivial  $\phi  \in 
 {L}^2[0,1]$  orthogonal to the range of  $DF(u,f)$, that is 
 
$$ 
\int_0^1 \phi \{(a \dot{u}_x)_x + f'(u)\dot{u} + \dot{f}(u) \} = 0 
\textrm{ }\forall \dot{u} \in {H}^2_N[0,1] \textrm{ and  } 
\forall \dot{f} \in {\cal C}^2(\Erre,\Erre). 
$$ 
Taking $\dot{f}=0$, 
$\int_0^1 \phi \{(a \dot{u}_x)_x + f'(u)\dot{u} \} = 0 
\textrm{ }\forall \dot{u} \in \mathnormal{H}^2_N(0,1)$ 
and   $\phi \in \mathnormal{H}^2_N(0,1) \cap \mathcal{C}^2(0,1)$ 
 is a weak, therefore, strong solution of (\ref{lve}). 
 
 Taking now, $\dot{u} = 0$ 
$\int_0^1 \dot{f}(u)\phi = 0 \textrm{ } \forall \dot{f} \in 
\mathcal{C}^2(\Erre,\Erre)$, as claimed.

 We are now in a position to prove our main result 
 {\em if additional properties on the equilibria are assumed}. 
 We first show how the results for {\em constant} coefficients  can be 
 obtained with our approach. 
 
  \begin{theo}  \label{hypconst} 
  Suppose the diffusion coefficient  in (\ref{parabolic}) is constant. Then 
 there exists a 
  residual set  $  { F}  $ in $  {\cal C}^2_S ( \Erre) $ such that, if 
   $f \in { F}  $, all equilibria of  (\ref{parabolic}) 
   are hyperbolic. 
  \end{theo} 
 
 \noindent{\bf  Proof.  } We prove that, for an open dense set of $f$ 
 the equilibria satisfying 
 $ \| u  \|_{H^2[0,1]} \leq N $ 
are all hyperbolic and then take intersection. 
 Again, openness  is no problem. 
 
Let    $Y_N$ be the open dense set  given by 
 (\ref{constequilibria}).

 We apply the Transversality Theorem for  the map 
 \begin{eqnarray} 
 \Psi:B_N \times Y_N & \mapsto & {L}^2[0,1] \label{Phi} \\ 
(u,f) & \mapsto  & (a u_x)_x + f(u) 
\end{eqnarray} 
 
  If $0$ is not a regular value of $\Psi$ it follows from lemma 
\ref{condition} that there is a solution  $u$ of (\ref{bvp}) 
 and a  nontrivial solution $\phi$ of the corresponding 
 linearized equation   (\ref{lve}) such that 
 $\int_0^1 \dot{f}(u)\phi = 0 \textrm{ } \forall \dot{f} \in 
\mathcal{C}^2(\Erre, \Erre)$. 
 
 But then, it follows from lemma  \ref{sumzero} that 
\[ 
 \sum_{p \in C_q \cap (0,1)} {\phi(p)}{ \sqrt{a(p)}} 
  + \frac{1}{2} \sum_{p \in C_q \cap \{0,1\}} 
  { \phi(p)} \sqrt{a(p)} = 0. 
 \] 
 
  Since $u_x$  and $\phi$ are two linearly independent  solutions 
  of the  second order  linear 
 differential equation in (\ref{lve}) their Wronskian must be a nonzero 
  constant, 
 that is, 
\[ u_x \phi - u_{xx} \phi = k \neq 0. 
 \] 
 At critical points $p$ of $u$,  we then have 
 $ u_{xx}(p) \phi(p) = k$ and also 
 $ u_{xx}(p) = - \frac{f(u(p))}{a}$. 
 and thus,  $\phi(p) = - \frac{k a}{f(u(p))}$. 
 
 But then, for any critical value $q $ of $u$,  the above sum gives 
  $ \frac{k  a^{3 \over 2}}{f(q)} = 0 $,  a contradiction.

 Therefore $0$ must be a regular point for $\Psi$. By the transversality 
 Theorem, it must be also a regular point for $\Psi(\cdot,f)$ for f 
 in a residual, therefore dense set of $ Y_N$ and thus, also of 
 ${\cal C}^2_S $.

 \begin{rem} 
  The argument above can be pushed a little further to cover the case 
 of monotonic $ ({\cal C}^1)$ coefficients. In fact, suppose $a'(x) \neq 0 $ 
 and let $q = f(0)$. If $p$ is another critical point of $u$, with 
 $u(p) = q$ then, as shown during the proof of lemma \ref{convex} 
 $ \frac{1}{2} \int_{0}^{p} a_x u_x^2 = 0.$ If $u$ is not constant, this 
 implies $ a_x = 0 $ in  $ [0,p] $. Therefore we can conclude, as above, 
 that $ \phi(p) = \phi(0)$ and obtain the same contradiction. 
 \end{rem} 
 
 For general diffusion coefficients, our arguments encounter difficulties for 
 solutions, which have  many  critical points  for each critical value. 
   We  first define them precisely.

 \begin{defin} \label{exceptional} 
 We say a nonconstant solution of (\ref{bvp}) is an exceptional 
  equilibrium  if 
\begin{enumerate} 
 \item  $u$ is not hyperbolic. \label{ex1} 
 \item    For any critical point $p$, there exists another 
 critical point $\bar{p}$ with the same value. \label{ex2} 
 \item There are  critical points $p,q$ with $u(p) = u(0)$, $u(q)=u(1)$, 
  such that $\phi(p) \neq \phi(0) $,   $\phi(q) \neq \phi(1) $, 
 for some (and therefore for any) nontrivial solution of  (\ref{lve}). 
 \label{ex3} 
 \end{enumerate} 
\end{defin} 
 We will also call {\em nonexceptional} any equilibrium that is not 
 exceptional. 
 
   For a given $f$ we   denote by 
 (\ref{bvp})     $E^{*} = E^{*}(a,f)$ the 
 set of exceptional solutions of (\ref{bvp}). 
 
   We can  prove that  $E^{*}$ 
   is really `exceptional' in the following sense:

 \begin{lemma} 
 For any $f \in {\cal C}^2_S ( \Erre) $ there exist at most a finite 
 number of exceptional equilibrium points $u$ with  $||u || \leq n$, where 
  $||u || $ 
 is the norm of $u$ in  $ H^{1}([0,1])$. 
 \end{lemma} 
 
 \noindent{\bf  Proof.  } 
  We  introduce    the linear variational equation around a solution 
 $u$ of (\ref{ivp}) 
 
 \begin{equation} 
  \left\{ 
    \begin{array}{lr} 
      {\displaystyle {(a (x) v_x)}_x 
     + f'(u)  v 
       \: = 0 }, \: &  0 \leq x \leq 1 \\ 
      v_x (0)= 0 \quad 
          v (0)= 1 & 
     \end{array} 
   \right.   \label{lvi} 
    \end{equation} 
 
  Then $ \frac{\partial u}{\partial u_0}$ the derivative of solutions 
 of (\ref{bvp}) with respect  to the initial condition is the 
 unique  solution of (\ref{lvi}).

  Suppose $u$ is a solution of  (\ref{bvp}). 
  and  $\phi$  a non-trivial 
  solution of (\ref{lve}).  Then 
 $ \frac{\partial u}{\partial u_0}$   of (\ref{lvi}) 
 around $u$ is a multiple of $\phi$. 
 Reciprocally, if 
  $ \frac{\partial u}{\partial u_0}$ 
 satisfies $  \left(\frac{\partial u}{\partial u_0}\right)_x(1)=0$ 
 then  $ \frac{\partial u}{\partial u_0}$  is  a  non-trivial  solution 
 of (\ref{lve}).

  Let    $u$ be an exceptional solution of (\ref{bvp}) and let 
 $ \phi \not\equiv 
  0  $ and $\frac{\partial u}{\partial u_0}$ the corresponding 
 solutions of ( \ref{lve}) and   (\ref{lvi}) respectively. 
 As observed above, we have $\frac{\partial u}{\partial u_0} = k \phi$, 
 where $k\neq 0 $ is a constant. 
 If $p \in [0,1]$ is such that $u(p) = u(0)$, $u'(p) = 0$, then 
 we write  $\bar{u}_0 = u(0)$  and denote by $ u(u_0, \cdot)$ the unique solution 
 of  ( \ref{ivp}) with  initial value $u_0$. Then, by the implicit function 
 theorem ($u$ is not constant!) there exists a unique 
 point $p(u_0)$ for $u_0 $ in a neighborhood of $\bar{u}_0$, such that 
 $ u_x(u_0,p(u_0) )= 0$. 
 We then have

 \begin{eqnarray*} 
  u(u_0, p(u_0)) & = & u(\bar{u}_0, p)+ 
\left( \frac{\partial u}{\partial u_0}(\bar{u}_0, p) 
  + u_x(\bar{u}_0, p)p'(u_0)\right)(u_0-\bar{u}_0)  \\ 
       & + & o(u_0-\bar{u}_0) \\ 
 & = &\bar{u}_0+  k \phi(p)(u_0-\bar{u}_0)+ o(u_0-\bar{u}_0) 
 \end{eqnarray*} 
 
 On the other hand 
 
\begin{eqnarray*} 
  u(u_0, 0) & = & u(\bar{u}_0, 0)+ 
\left( \frac{\partial u}{\partial u_0}(\bar{u}_0, 0) 
  \right)(u_0-\bar{u}_0) + o(u_0-\bar{u}_0) \\ 
 & = & \bar{u}_0 + k \phi(0)(u_0-\bar{u}_0)+ o(u_0-\bar{u}_0) 
 \end{eqnarray*}

 Since, by hypothesis $ \phi(p) \neq \phi(0)$ 
 we have,  for any solution of (\ref{ivp}) with initial condition $u_0$ close 
 to $ \bar{u}_0$ 
 $ u(u_0, p(u_0)) \neq  u(u_0, 0)$.

 Of course, $ u(u_0, \cdot)$ 
 may be a solution of (\ref{bvp}) but it will not   be an exceptional 
 solution. 
 Thus, such solutions must be in a discrete set. Since the set of 
 solutions of (\ref{bvp}  ) with $ ||u|| \leq n$ is also compact 
  in $H^1[0,1]$., 
 and the set of exceptional solutions is closed in   $H^2[0,1]$ 
  it must be finite as claimed.

   Assuming analyticity     we can show  generic hyperbolicity for 
 the  nonexceptional equilibria. 
 
  \begin{lemma} \label{nonexcep} 
  Suppose the diffusion coefficient $a$  and $f$ in (\ref{parabolic}) is 
 analytic 
 and it is not an even function about the point $x = \frac{1}{2}$ 
   Then, there is a 
  residual set  $  { F}  $ in $  {\cal C}^2_S ( \Erre) $ such that, if 
   $f \in { F}  $, all nonexceptional  equilibria of  (\ref{parabolic}) 
   are hyperbolic. 
  \end{lemma}

  \noindent{\bf  Proof.  }  We proceed as in theorem  \ref{hypconst} but 
 now restricting $u$ to the (open) subset of functions $u$ in  $B_N$ 
 which do not satisfy at least one of the conditions in definition 
 (\ref{exceptional})  (that it is, excluding the potential 
 nonexceptional equilibria).  Suppose $u$ is a critical point of 
 (\ref{Psi}), and $ \phi$ is a solution of  (\ref{lve}), 
 with $ \int_0^1 f(u) \phi = 0 $ for any continuous $f$. 
  If $u$ fails to satisfy (\ref{ex3}) of definition 
 \ref{exceptional} then, it follows from lemma \ref{sumzero} 
 that $\phi(0) = 0 $, so $\phi \equiv 0$, a contradiction. 
 Suppose then that $u$ does not satisfy (\ref{exceptional}-\ref{ex2}) 
   and let $p$ be the 
  required critical point of $u$, with 
 no other critical point at the same level. As in the proof of 
 lemma \ref{orbit} we denote by 
  $ 0 = x_0 < x_1 <  x_2 < \cdots < x_n = 1$ the critical points of $u$ 
 and   by  $\xi_i$, $ i=1,2 \cdots, n$ be the inverse of 
 $u$ in the interval $[x_{i-1},x_{i} ]$. Suppose that  $p = x_l$ is a local minimum 
(for local maximum, the argument is similar).  If $q$ is a regular value of $u$ 
   and $ \xi_{i_1},  \xi_{i_2}, \cdots, \xi_{i_n}$  are its inverse images under $u$, 
 and $ \Phi = \frac{\phi}{u_x}$ we 
 have, by lemma \ref{sumzero} 
  $ S(q) = 
 \sum_{j=1}^{n} \Phi(\xi_{i_j})(q) \mbox{\rm sign} {u'(\xi_{i_j})}(q) = 0$. 
 
 When $q$ passes through $u(p)$, coming  from below, the term 
 $\Phi(\xi_{l})(q) \mbox{\rm sign} {u'(\xi_{l})}(q) + 
    \Phi(\xi_{l+1})(q) \mbox{\rm sign} {u'(\xi_{l+1})}(q) $ 
 is added to 
  the above sum. Therefore, we must have 
  $\Phi(\xi_{l})(q) \mbox{\rm sign} {u'(\xi_{l})}(q) + 
    \Phi(\xi_{l+1})(q) \mbox{\rm sign} {u'(\xi_{l+1})}(q) = 0 $, for 
 $q$ slightly above $u(p)$. 
   By analyticity, this  should also be  true for $  u(p) <q < u(x_{l-1})$, if 
$ u(x_{l-1}) \leq  u(x_{l+1})$ or  $  u(p) <q < u(x_{l+1})$ otherwise. 
 Suppose wolog, we are in the first case. Let  $a = (x_{l-1})$, and $b$ the 
 unique point  in $ \left]  (x_{l-1}), (x_{l+1}) \right[$ with $u(a) = u(b)$, 
 From( the second part of) lemma \ref{sumzero}, we have then 
   $ \int_a^b f(u) \phi = 0 $ for any continuous $f$.

  From lemma \ref{inorbit}, it follows that   $a = g(u) $  for some analytic 
 function $u$ in the interval $ [a,b] $ and, in particular 
  $a_x  = 0$ at the  critical point $p$. 
(this can also be obtained 
 from  $a_x u_x \in {\cal O}(u)$ or  $a_x \phi \in {\cal O}(u)$).

  Write  $\xi_1:  [u(p), u(a)] \mapsto [a,p] $
         $\xi_2:  [u(p), u(b)] \mapsto [p,b]$ 
 for  the inverse of 
 $u$ in each interval. From  (\ref{16}) of lemma \ref{orbit}, it follows also that 
 $u_x(\xi_1(y)) = - u_x(\xi_2(y))$ for any $ y \in [u(p), u(a)]$, hence 
  $ \xi_1'(y)) = - \xi_2'(y)$ for any $ y \in [u(p), u(a)]$. If  $\xi_1(y) = x_1$, 
  $x_2(y) = x_2$, then 
  $ x_1 - p =   \xi_1(y) - \xi_1(u(p)) = \int_{u(p)}^{y} \xi_1'(y) \, dy = 
   -  \int_{u(p)}^{y} \xi_2'(y) \, dy = - ( x_2 - p ) $. 
  Thus $ p - \xi_1(y) =  \xi_2(y) - p $    for any $ y \in [u(p), u(0)]$. 
  Therefore, $a$ is an even function about the point $p$ in the interval 
 $ [a,b]$. But, being analytic, $a$ must be an even function 
 about  $p$ in  the whole interval $[0,1]$. In particular $p = {1 \over 2}$, 
 contradicting the hypothesis.

 \vspace{3mm} 
 
   This result is not completely satisfactory,  since it imposes 
 conditions on the (unknown) solutions of the problem. The natural question 
 is then if one could rule out beforehand the existence of nonexceptional 
 solutions. It turns out that this is possible, if $a$ has few oscillations. 
 
  We say that a ${\cal C^1}$ function $a : [0,1] \mapsto \Erre $ has $n$ intervals of 
 monotonicity if there exists a partition $ 0 = x_0 < x_1 < \cdots < x_n =1$, 
 such that $a'(x)$ is strictly positive or negative in each subinterval, 
 $[x_{i-1}, x_{i}]  $.

  \begin{lemma} \label{convex} 
 Suppose $a$ has at most two  intervals of monotonicity  in $[0,1]$. Then 
   there are no exceptional    solutions of (\ref{bvp}). 
   \end{lemma} 
 
 \noindent {\bf Proof. } 
 If $p, \bar{p} \in  [0,1]$ are two critical points with the same value, then 
 multiplying (\ref{bvp}) by  $u_x$ and integrating, we get 
 
   \begin{eqnarray}\label{intzero} 
   0 & = &  \int_{p}^{\bar{p}}  \{(au_x)_x\} u_x + f(u) u_x \,  d \, x \nonumber \\ 
     & = &  ({au_x}^2(\bar{p}) -  ({au_x}^2(p) - 
 \frac{1}{2} \int_{p}^{\bar{p}} a \frac {d}{dx}  {(u_x)}^2   \, dx \nonumber  \\ 
     & = &   -  \frac{1}{2}  {(au_x)}^2(\bar{p})   +  \frac{1}{2}  {(au_x)}^2(p) 
  +  \frac{1}{2} \int_{p}^{\bar{p}} a_x  {(u_x)}^2   \, dx \nonumber  \\ 
   & = &    \frac{1}{2} \int_{p}^{\bar{p}_2} a_x  {(u_x)}^2   \, dx 
     \end{eqnarray}

 By hypothesis  we cannot 
 have a critical point $p$ of $u$ in 
  $ \{0,1\}$ with $C_{u(p)} = \{p\}$ (that is, 
 such that no other critical point with the same value exists).

 Suppose   there are three critical points  $p_1< p_2 < p_3$, with the same value $q$, 
 Then, if  $ p_2 \leq c$ 
    $a_x$ has constant sign in $[p_1,p_2]$  so 
 $ \int_{p_1}^{p_2}   a_x  {(u_x)}^2   \,  d \, x \neq 0 $, contradicting 
 (\ref{intzero}). 
   By a similar argument, we also  cannot have   $ p_2 \geq c$, so there are 
 at most  and therefore, exactly two points in   $u^{-1}(q)$. In 
 particular, the number of critical points must be even

  If $u(0) \neq u(1)$, we  must have  two critical points     $p_0$ and  $p_1$ 
in $(0,1)$ with $u(p_0) = u(0) $, $u(p_1) = u(1) $. 
 If $ p_0 < c$ or $p_1 > c $  then, since  $a_x$ has constant sign in at least 
  one of the intervals $ [0,p_0],  [p_1,1] $ then 
 $  \int_{0}^{p_0} a_x  {(u_x)}^2   \, dx \neq  0$ or 
  $  \int_{p_1}^{1} a_x  {(u_x)}^2   \, dx \neq 0$, in contradiction with 
 (\ref{intzero}). 
 
 If $p_0 > c > p_1 $ then, again by (\ref{intzero}) 
  $ 0 =  \int_{0}^{p_0} a_x  {(u_x)}^2   \, dx =  \int_{0}^{p_1} a_x  {(u_x)}^2   \, dx 
 +  \int_{p_1}^{p_0} a_x  {(u_x)}^2   \, dx$, and 
   $ 0 =  \int_{p_1}^{1} a_x  {(u_x)}^2   \, dx = 
  \int_{p_1}^{p_0} a_x  {(u_x)}^2   \, dx  +  \int_{p_0}^{1} a_x  {(u_x)}^2   \, dx$. 
 Subtracting, we obtain 
   $ 0 =  \int_{p_0}^{1} a_x  {(u_x)}^2   \, dx - 
  \int_{0}^{p_1} a_x  {(u_x)}^2   \, dx = 0 $. 
 But this is not possible, since $a_x$ has different signs in the intervals 
 $  [0,p_1]  $, and     $  [p_0, 1]  $. 
 
 If  $u(0) = u(1)$, since the zeroes of  $u$ are nondegenerate, we must 
 have an odd number of critical points, contradicting the 
 conclusion above. This proves the result.

 As an immediate consequence, we have

  \begin{cor} 
  Suppose the diffusion coefficient $a$  and $f$ in (\ref{parabolic}) is analytic, 
  and has only two intervals of monotonicity in $[0,1]$, and is not even about 
 the point $x = {1 \over 2}$. Then 
 there exists a 
  residual set  $  { F}  $ in $  {\cal C}^2_S ( \Erre) $ such that, if 
   $f \in { F}  $, all equilibria of  (\ref{parabolic}) 
   are hyperbolic. 
  \end{cor}

 \noindent {\bf Acknowledgment: \ } Part of this work was carried out while the 
 author was visiting the Instituto Superior T\'ecnico de Lisboa, Portugal. 
 He  wishes to acknowledge   the support and  warm hospitality of the 
 Institution. He also wishes to thank 
 Prof. Carlos Rocha and M. C. Pereira for  many  discussions and suggestions.

 \end{document}